\documentclass{article}
\usepackage{geometry}[margin=1in]
\usepackage[utf8]{inputenc}
\usepackage[maxbibnames=99, backend=biber,style=alphabetic]{biblatex}

\usepackage{comment}
\usepackage{amsmath}
\usepackage{amssymb} 
\usepackage{amsfonts}
\usepackage{amsthm}
\usepackage{bbm}
\usepackage{mathtools}
\usepackage{thmtools, thm-restate}
\usepackage{hyperref}
\usepackage{tikz}
\usepackage{tikz-cd}
\usetikzlibrary{calc}
\usetikzlibrary{positioning}
\usepackage{pgfplots}

\usepackage{algorithm}
\usepackage{algpseudocode}
\usepackage{subfig}
\usepackage{wrapfig}

\newtheorem{theorem}{Theorem}

\newtheorem{lemma}{Lemma}

\newtheorem{definition}{Definition}

\newtheorem{claim}{Claim}
\newtheorem{assumption}{Assumption}

\newcommand{\R}{\mathbb{R}}

\newcommand{\E}{\mathbb{E}}

\DeclareMathOperator{\vol}{vol}
\DeclareMathOperator{\rank}{rank}

\pgfplotsset{compat=1.18}
\bibliography{references}
\title{Convergence to Radial Symmetry in Iterative Convolution-Thresholding Dynamics
}
\author{
Mirabel Reid\thanks{School of Computer Science, Georgia Institute of Technology (mirabel.reid@tu-darmstadt.de, dzhang381@gatech.edu)}
\and 
Daniel J. Zhang\footnotemark[1]
}
\date{}
\begin{document}
\maketitle
\begin{abstract}
We study a discrete-time spatially extended dynamical system motivated by the continuum limit of binary neuron networks on geometric random graphs. The model evolves a function $\psi:\R^d\rightarrow[0,1]$ by iterative smoothing and sharpening; that is, $\psi_{t+1} = \gamma_t\circ (g * \psi_t)$, where $g$ is a radial convolution kernel and $\gamma_t$ is a monotone map. Under mild regularity conditions on $g$ and $\gamma_t$, we prove that $\psi_t$ tends toward radial symmetry as $t \rightarrow \infty$. Notably, a special case of this process recovers the Merriman–Bence–Osher (MBO) scheme for the motion of interfaces by mean curvature, and we provide a novel analysis of its behavior. Our results connect the dynamics of spatial binary neuron networks with classical models of interface motion, and we establish general conditions under which spatial dependence drives activity toward radial symmetry.
\end{abstract}

\section{Introduction}
Spatially extended dynamical systems, in which state and interaction are dependent on location in physical space, arise in modeling across biology, physics, materials science, and other physical sciences. Many such spatially extended systems are governed by interactions that are both local and symmetric: particles influence each other most strongly at short distances, and the strength of the influence depends only on their relative separation. Such isotropic interactions arise naturally in models of diffusion, aggregation, and spatially embedded networks, and a central challenge in these systems is to understand how these local interactions lead to macroscopic organization. 

A fundamental discrete-time model capturing this behavior is \textit{iterative convolution-thresholding dynamics}. In this framework, a binary state variable is evolved in two phases: local influence is computed by convolving the state with a fixed kernel, and the state is updated by taking a threshold. The result is a gradual `smoothing' effect, and it has been applied to image processing and segmentation~\cite{merkurjev2013mbo, wang2022iterative} as well as simulating interface motion~\cite{merriman1992diffusion}.

%\dz{This makes it sound like our main contribution is analyzing a generalization of the model. If I remember correctly, we show something new for the thresholding model too, which was the ``main'' result in the sense that the reader is more likely to care. (Unless you have an interesting example of a monotone map that is not thresholding.}
Our main contribution is a formal analysis of the evolution of the state variable as $t\rightarrow \infty$. Specifically, we study a natural generalization of this model in which the thresholding step is replaced by a monotone nonlinearity. The state at $t$ is a function $\psi_t:\R^d\rightarrow[0,1]$, and the evolution is given as

\begin{equation}
\psi_{t+1} = \gamma_t\circ (g * \psi_t)
\label{eq:dynamics}
\end{equation}

where $g$ is an isotropic convolution kernel and $\gamma_t$ is a (possibly time-dependent) monotone map. This formulation interpolates between the sharp interfaces elicited by thresholding and the diffuse activation functions common in biological applications. 

This class of dynamical systems arises in a number of physical applications. When $\gamma_t$ is the threshold at 1/2, the resulting dynamics correspond to the Merriman-Bence-Osher (MBO) scheme to simulate the motion of an interface by mean curvature~\cite{merriman1992diffusion}. Typically $g$ is the heat kernel with width $\tau$, $g(x) = G_\tau(x) = (4\pi \tau)^{-d/2}\exp\left(-|x|^2/4\tau\right)$, and in this case, the MBO scheme has been shown to converge to mean curvature flow as $\tau \rightarrow 0$\cite{evans1993convergence, swartz2017convergence}. 

The process described by Equation \ref{eq:dynamics} is also motivated by the continuum limit of a geometric network of binary neurons. Since neural connections are governed by spatial and material constraints, spatial models such as the geometric random graph are sometimes used to model the brain's graph structure~\cite{duchemin2023random,bullmore2009complex, barthelemy2011spatial}. In the limit as the number of neurons approaches infinity, Equation \ref{eq:dynamics} can approximate the evolution of firing probabilities in a discrete-time, binary neuron model on a geometric random graph with edge probabilities determined by $g$. Reid and Vempala~\cite{reid2023k} used this continuum limit to study the $k$-cap process, a model in which the neurons which fire at each time step are determined by a $k$-winners-take-all function, taking $\gamma_t$ as a volume-preserving threshold. However, the discussion of the model was limited to $d=1$, and deeper analysis was left to future work. 

In this work, we study a general class of isotropic convolution-thresholding dynamics, and we utilize a novel application of the method of moving planes~\cite{alexandrov1962uniqueness} to characterize the limiting structure. Our main results are as follows. (i) Under mild regularity conditions on $g$ and $\gamma_t$, we prove that $\psi_t$ tends toward radial symmetry as $t \rightarrow \infty$. (ii) In the special case where $g$ is analytic and $\gamma_t$ is a volume-preserving threshold, we prove that $\psi_t$ converges to the indicator function of a ball in $\R^d$, resolving the open question of~\cite{reid2023k}. In particular, this establishes that the volume-preserving threshold step, previously studied in the context of the MBO scheme for volume-preserving motion by mean curvature~\cite{ruuth2003simple,laux2017convergence}, is well defined for analytic kernels. 

\subsection{Definition}
Let $\psi_t:\R^d \rightarrow [0,1]$ be an $L_1$-integrable function representing a scalar measure of activity at a point $x$ at time $t$. The evolution of $\psi_t(x)$ is based on the activity of its neighborhood, with nearer points exerting more influence. Formally, the update $\psi_{t+1}$ is defined by \textit{convolving} $\psi_t$ with a kernel $g$, followed by \textit{composing} with a monotone function $\gamma_t$ (Definition~\ref{alg:cts-general}).
\begin{definition}[Generalized Iterative Convolution-Thresholding Process]
\label{alg:cts-general}
\leavevmode\\
Let $\psi_0 \in L^1(\mathbb{R}^d)$ with values in $[0,1]$. We define the interaction field
\[
f_t := \psi_t * g,
\]
and evolve $\psi_t$ according to
\[
\psi_{t+1} = \gamma_t\circ f_t, \quad t \ge 0.
\]
where $g : \mathbb{R}^d \to [0,1]$ is a kernel and $\gamma_t : \mathbb{R}^+ \to [0,1]$ is a sequence of monotone increasing functions.
\end{definition}
We make the following assumptions on the sharpening function $\gamma_t$. In particular, we will assume that $\gamma_t$ is defined such that the support of $\psi_t$ is bounded. 
\begin{assumption}
\label{assumption:gamma}
Assume that $\gamma_t:\R^+\rightarrow [0, 1]$ is monotone increasing, and satisfies $\gamma_t(f_t(x)) = 0$ whenever $\|x\|>R$.
\end{assumption}
\noindent We also make the following regularity assumptions on the kernel function $g$.

\begin{assumption}
\label{kernel-assumptions}
Assume that the kernel function $g:\R^d\rightarrow[0,1]$ has the following properties. 
\begin{itemize}
\item $g(x) = \tilde{g}(\lVert x \rVert)$ for a $\tilde{g}:\R_{\ge 0}\rightarrow [0,1]$ (i.e. the function is radially symmetric).
\item $g(0) = 1$ and $g(x)> 0$ for all $x\in \R^d$.
\item $g$ is $L-$Lipschitz continuous and has continuous first and second derivatives. 
\item $\tilde{g}'(x) < 0$  for all $x>0$. 
\end{itemize}
\end{assumption}

We can assume without loss of generality that $\lim_{\| x \| \rightarrow \infty} g(x) = 0$, by the following argument. The other assumptions together imply that $\lim_{\| x \| \rightarrow \infty} g(x) = l\ge0$ (i.e., the limit exists). If $l \neq 0$, letting $\hat{g}(x-y) = g(x-y) - l$, $f_t(x) = \int_{\R^d}\psi(y)g(x-y)\, dy = \int_{\R^d} \psi(y)\hat{g}(x- y)\, dy - l\int_{\R^d}\psi(y)\,dy$. Shifting $f_t$ by a constant does not affect the dynamics, as long as $\gamma_t$ is shifted in kind.

This definition encompasses many commonly used kernel functions. As a motivating example, let $g$ be the Gaussian function $g(x) = (2\pi\sigma^2)^{-d/2}\exp\left(-\lVert x \rVert^2\right/2\sigma^2)$. Another common example is an inverse-square distance function, $g(x) = \frac{1}{1+\| x\|^2}$.

This problem is motivated by the following special case; suppose that $\psi_t$ is the indicator function of a compact set $A_t \subset \R^d$ whose volume $|A_t| = \alpha$ is preserved over time. This can be maintained by choosing $\gamma_t$ to be a time-varying threshold, $\gamma_t(x) = \mathbbm{1}\{x \ge C_t\}$. The threshold $C_t$ is then chosen such that $|A_t| = |\{x: f_t(x) \ge C_t\}| = \alpha$ (we will show in Lemma \ref{lem:analytic} that, under the additional condition that the kernel $g$ is analytic, such a threshold always exists). Since the set $A_t$ is chosen by picking the `top' $\alpha$ set according to the interaction field $f_t$, we call this the \textit{$\alpha$-cap process}.

\begin{definition}[$\alpha$-cap Process]
\label{alg:cts}
Let $A_0 \subset \R^d$ be a compact set. We define the interaction field
\[
f_t := \psi_t * g,
\]
where $g : \mathbb{R}^d \to [0,1]$ is a kernel. Define the volume-preserving threshold $C_t = \inf\{ C\in \mathbb{R}^+ : |\{x : f_t(x) \ge C\}| < \alpha\}$. Then evolve $A_t$
\[
A_t = \{x : f_t(x) \ge C_t\}
\]
\end{definition}

\begin{wrapfigure}{R}{0.3\textwidth}
\centering
\includegraphics[width=0.29\textwidth]{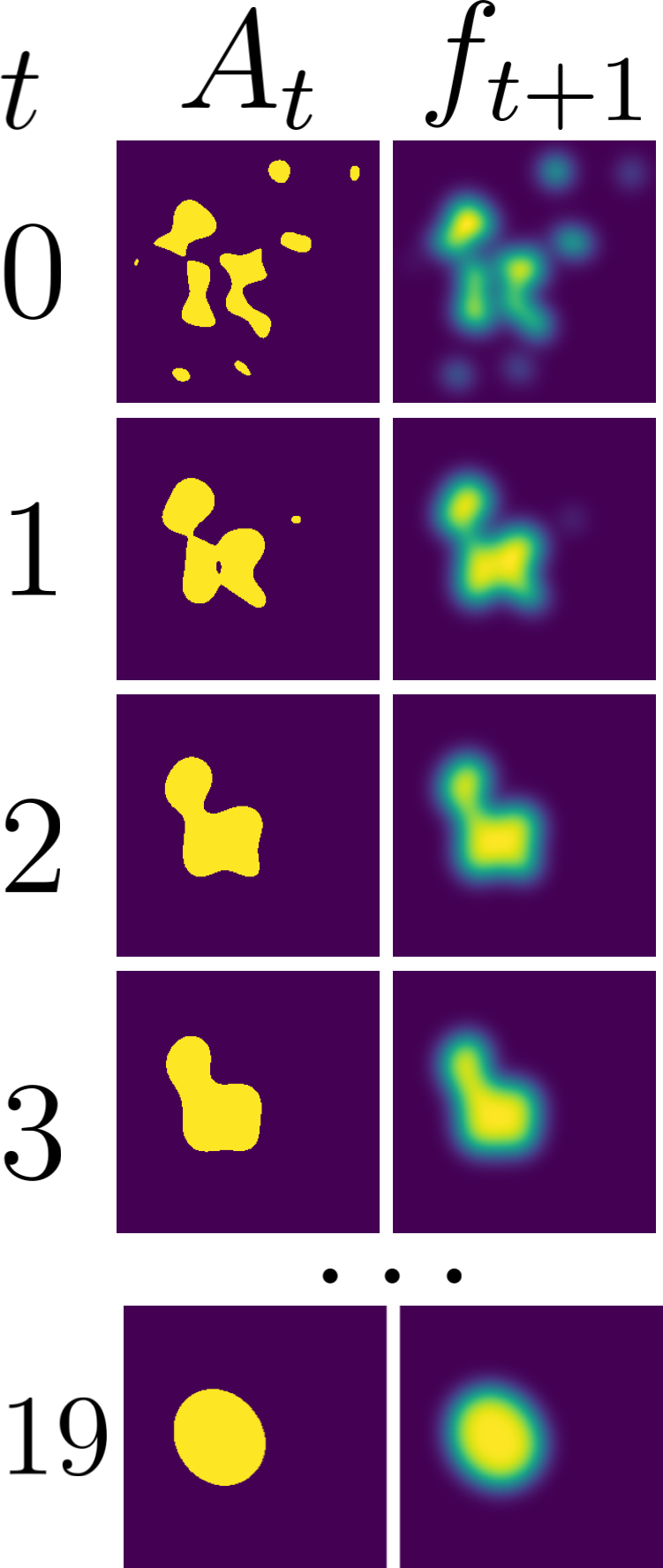}
\caption{\label{fig:example1}The evolution of a set in $\R^2$ by the $\alpha$-cap process (Definition \ref{alg:cts}), with $g$ set to the Gaussian kernel. }
\end{wrapfigure}

An example of the evolution of a set by the $\alpha$-cap process is given in Figure \ref{fig:example1}. Although the initial set is disconnected and non-convex, it is gradually smoothed, until converging to a circle. Indeed, as we will prove in Section \ref{sec:volume-preserving}, for a sufficiently regular kernel $g$, all initial sets $A_0$ converge to a ball in $\R^d$ as $t\rightarrow \infty$. 

\subsubsection{Relationship to Binary Neuron Communication}
The $\alpha$-cap process (Definition \ref{alg:cts}) can be viewed as a continuous approximation for binary neuron communication on graphs with geometric structure. This is motivated by the $k$-cap process, a model for neuron communication in a setting where the global firing rate at each time step is fixed~\cite{reid2023k}. At each time step $t>0$, each vertex is assigned a state {\em active} (1) or {\em inactive} (0). The {\em active set} $A_t\subseteq V$ consists of the $k$ vertices with the highest degree in $A_{t-1}$ (with ties broken randomly). 

\begin{definition}[$k$-cap Process]
\label{def:kcap}
Let $G$ be a graph on $n$ vertices, and let $A_0\subset [n]$ with $|A_0| = k$ be the set of initial active vertices.

We define the interaction field
\[
f_G(x; A_t) :=  \frac{1}{n-1}\sum_{i\in [n]\setminus\{x\}} \mathbbm{1}\{i\in A_t\}e_{G}(i, x)
\]
Then, for all $t\ge 0$, update $A_t$ according to 
\[A_{t+1} = \texttt{topk}(\{f_G(x; A_t)\}_{x\in [n]})\]
with ties broken randomly. 
\end{definition}

Motivated by the spatial structure of biological neural networks, \cite{reid2023k} studied this process on Gaussian geometric random graphs (Definition \ref{def:geometric}, with $g$ being the Gaussian kernel). It was found that, under a specific parameterization, the active sets $A_t$ converge to lie within a small ball.

\begin{definition}[Soft Geometric Random Graph]
\label{def:geometric}
Let $\mathcal{X} \subset \mathbb{R}^d$ be a compact set, and let $g: \mathbb{R}_{\geq 0} \rightarrow [0,1]$. A \textit{soft geometric random graph} $G=G(n, \mathcal{X}, g)$ is defined as follows:
\begin{itemize}
\item Let $V\subset\mathcal{X}$ be a set of $n$ vertices chosen uniformly and independently at random from $\mathcal{X}$
\item For each pair $x,y\in V$, $x \neq y$, let $e_{G}(x, y)=1$ independently with probability $g(\|x-y\|_2)$. 
\end{itemize}
\end{definition}

We can derive the interaction field $f$ of the $\alpha$-cap Process (as given in Definition~\ref{alg:cts}) by taking the limit of $f_{G_n}$ (as given in Definition~\ref{def:kcap}) on a sequence of graphs $\{G_n\}$ drawn from the soft geometric random graph model. We abuse notation and define $g:\mathcal{X}\times\mathcal{X}\to[0,1]$ by $g(x,y)=g(x-y)$.

Since the domain of $f_{G_n}$ varies with $n$, the limit must be defined carefully. We prove the convergence of $\{f_{G_n}\}$ to the interaction field of Definition~\ref{alg:cts} via the following theorem:
\begin{theorem}
    Let $A\subseteq\mathcal{X}$. For $n\ge 2$, let $G_n\sim G(n,\mathcal{X},g)$ independently, and $V_n=V(G_n)$.

    Define
    \begin{align*}
        f_n(x)&=\dfrac{1}{n-1}\sum_{y\in V_n\setminus\{x\}}\mathbbm{1}_{A}(y)e_{G_n}(x,y)\\
    f_\infty(x)&=\int_\mathcal{X}\mathbbm{1}_{A}(y)g(x,y)dy
    \end{align*}
    
    Then, for any measurable $\phi:\mathcal{X}\to[0,1]$,
    \begin{align*}
        \dfrac{1}{n}\sum_{x\in V_n}\phi(x)f_n(x)&\to \int_{\mathcal{X}}\phi(x)f_{\infty}(x)dx\text{ almost surely}
    \end{align*}
    \label{thm:cts-limit}
\end{theorem}

\begin{proof}

Fix $\phi:\mathcal{X}\to[0,1]$. Let $Y_n=\dfrac{1}{n}\sum_{x\in V_n}\phi(x)f_n(x)$ and $Z=\int_{\mathcal{X}}\int_{\mathcal{X}}\phi(x)\mathbbm{1}_A(y)g(x,y)dxdy$.

\begin{claim}
    $\E[Y_n]=Z$ and $\Pr(|Y_n-Z|\ge \epsilon)\le 2\exp(-\dfrac{\epsilon^2n}{4})$ for all $n\ge 2$ and $\epsilon>0$.
\end{claim}
\begin{proof}[Proof of claim]
We can construct $G_n\sim G(n,\mathcal{X},g)$ the following way:

Choose $x_1,\ldots,x_n\in \mathcal{X}$ and $t_{ij}\in[0,1]$ for $1\le i<j\le n$ uniformly and independently. Note that the $x_i$ are distinct almost surely. Let $t_{ji}=t_{ij}$.
Let $V_n=\{x_1,\ldots,x_n\}$ and $e_{G_n}(x_i,x_j)=1[t_{ij}\le g(x_i,x_j)]$. Now, define the following random variables $Y_n$.
\begin{align*}
    f_n(x_i)&=\dfrac{1}{n-1}\sum_{j\in[n]\setminus\{i\}} \mathbbm{1}_A(x_j)1[t_{ij}\le g(x_i,x_j)]\\
    Y_n&=\dfrac{1}{n}\sum_{i\in[n]}\phi(x_i)f_n(x_i)=\dfrac{1}{n(n-1)}\sum_{\substack{i,j\in[n]\\i\ne j}}\phi(x_i)\mathbbm{1}_A(x_j)1[t_{ij}\le g(x_i,x_j)]
\end{align*}
Note that the expectation of $Y_n$ is given by:
\begin{align*}
    \E[Y_n]&=\dfrac{1}{n(n-1)}\sum_{\substack{i,j\in[n]\\i\ne j}}\E\left[\phi(x_i)\mathbbm{1}_A(x_j)1[t_{ij}\le g(x_i,x_j)]\right]\\
    &=\int_{\mathcal{X}}\int_{\mathcal{X}}\phi(x)\mathbbm{1}_A(y)g(x,y)dxdy=Z
\end{align*}
If an individual $t_{ij}=t_{ji}$ is modified, it can only change $Y_n$ by at most $\dfrac{2}{n(n-1)}$. Similarly, changing $x_i$ can only change $Y_n$ by at most $\dfrac{2}{n}$. Hence, by McDiarmid's inequality \cite{mcdiarmid1989method},
\begin{align*}
    \Pr(|Y_n-\E[Y_n]|\ge\epsilon)&\le 2\exp(-\dfrac{2\epsilon^2}{n\cdot (\frac{2}{n})^2+\frac{n(n-1)}{2}\cdot(\frac{2}{n(n-1)})^2})\\
    &\le 2\exp(-\dfrac{2\epsilon^2}{8/n})= 2\exp(-\dfrac{\epsilon^2n}{4})
\end{align*}
\end{proof}

We now use the claim to show that $Y_n\to Z$ almost surely. We have
\begin{align*}
    \sum_{n=2}^{\infty}\Pr(|Y_n-Z|\ge n^{-1/3})\le 
    \sum_{n=2}^{\infty}2\exp(-\dfrac{n^{1/3}}{4})<\infty
\end{align*}
since for large enough $n$, $\exp(-\dfrac{n^{1/3}}{4})\le \exp(-2\log n)=\dfrac{1}{n^2}$, and $\sum_{n=1}^{\infty}\dfrac{1}{n^2}<\infty$.

Thus, by the first Borel-Cantelli lemma (see e.g. \cite{billingsley2017probability}), the probability that $|Y_n-Z|\ge n^{-1/3}$ infinitely often is zero.

Hence, $Y_n\to Z$ almost surely.\\
\end{proof}

Therefore, the interaction field of the $\alpha$-cap process (convolution of the set $A_0$ with the kernel function $g$) arises from the limit of the interaction field on finite geometric random graphs. 

Similarly, the sharpening function $\gamma_t$ models the firing probability of a vertex, as a function of the input from the previous active set. In the fixed-volume $\alpha$-cap process, $\gamma_t(x) = \mathbbm1\{x\ge C_t\}$ serves as an analogous function to the $\texttt{topk}$ function of the $k$-cap process (under the assumption that $|A_n|/n\rightarrow \alpha$). In general, by varying $\gamma_t$, the generalized convolution-thresholding dynamics can approximate different methods of determining firing probabilities as a function of input. For example, by scaling $\gamma_t$, the global firing rate can be held constant or allowed to vary over time. 
\subsection{Notation}
Unless otherwise stated, we operate in $\R^d$. Let $\mathbb{S}^{d-1}$ denote the set of unit vectors in $\R^d$. Let $\lVert x \rVert = \lVert x\rVert_2$ be the Euclidean norm of $x$.

We use $B(a,r)=\{x\in\R^d:\|x-a\|<r\}$ to denote open ball of radius $r$ around $a$.

We use $\mathbbm{1}_A$ to denote the indicator function for a set $A$, i.e. $\mathbbm{1}_A(x)=1$ for $x\in A$ and $\mathbbm{1}_A(x)=0$ for $x\notin A$.

To prove the convergence of the process (Definition~\ref{alg:cts-general}), we will employ a geometric argument, taking advantage of the symmetries of the kernel $g$. In order to facilitate this, we define convenient shorthand notation for operations with hyperplanes. 
\begin{definition}
    A hyperplane is
    \begin{align*}
        H=\{x:\langle x,u\rangle=a\}\subset \R^d 
    \end{align*}
    where $u\in \mathbb{S}^{d-1}$ and $a\in\R$. Let $H^-=\{x:\langle x,u\rangle<a\}$, and $H^+=\{x:\langle x,u\rangle>a\}$.
\end{definition}
Note that in this definition, hyperplanes are oriented. In other words, $H_1=\{x:\langle x,u\rangle=a\}$ is treated distinctly from $H_2=\{x:\langle x,-u\rangle=-a\}$, although the sets are the same.

For a point $x\in \R^d$, let $x^H$ denote the reflection of $x$ across the hyperplane $H$. If $H = \{y : \langle y, u\rangle = a\}$, the reflection has an explicit formula: $x^H = x-2(\langle x, u\rangle - a)u$. For a set $A \subset \R^d$, $A^H =  \{x^H : x\in A\}$ is the reflection of $A$ across $H$. When $H$ is clear from context, we may denote $A^+ = A \cap H^+$ and $A^- = A \cap H^-$. Similarly, define the function $\psi^H$ by $\psi^H(x) = \psi(x^H)$. 

\subsection{Related Work}
\paragraph{The MBO Scheme}A notable special case of the dynamics we study (Equation \ref{eq:dynamics}) arises from diffusion-generated motion by mean curvature. \textit{Motion by mean curvature} (MMC), also called \textit{mean curvature flow}, is a type of geometric flow that describes a surface which moves in the direction of its normal vector with velocity equal to its mean curvature. Merriman, Bence, and Osher~\cite{merriman1992diffusion, merriman1994motion} proposed a simple approximation scheme for MMC which evolves the surface by alternately diffusing and sharpening. Given an initial set $\Omega_0$ with characteristic function $\chi_0$, the surface is evolved in two steps:
\begin{align*}
&\text{Linear Diffusion }\:&\chi' &= \chi_t * G_\tau\\
&\text{Thresholding} \:&\chi_{t+1} &= \{\chi' \ge 1/2\}
\end{align*}
where $G_\tau$ is the Gaussian kernel parameterized by width $\tau$. Subsequent theoretical analysis has verified the convergence of the scheme to MMC as $\tau \rightarrow 0$~\cite{evans1993convergence, barles1995simple, swartz2017convergence}. Of particular relevance to our work, Ruuth and Wetton~\cite{ruuth2003simple} proposed a variant scheme which uses an adaptive threshold to approximate \textit{volume-preserving} motion by mean curvature. The convergence of this scheme in the zero-width limit was proven by Laux and Swarz~\cite{laux2017convergence}. Discrete variants of the MBO scheme have been applied to problems such as image classification~\cite{wang2022iterative,merkurjev2013mbo}, community detection~\cite{hu2013method}, and data clustering~\cite{laux2023large,merkurjev2014diffuse}.

\paragraph{Binary Neuron Models}As discussed in the previous section, a parallel application of this category of dynamics is in modeling inter-neuron communication. Discrete-time binary neuron models are a foundational abstraction in neuroscience, with roots going back to the McCulloch-Pitts neuron~\cite{mcculloch1943logical}. In such models, neurons are assigned a state --- either active (1) or inactive (0)--- and updated depending on the activity of neighboring neurons, typically via a threshold rule. The Hopfield model~\cite{hopfield1982neural}, an influential model of associative memory, couples binary neurons with Hebbian synaptic weights, admitting an energy function whose local minima correspond to stored patterns. It has been generalized by varying the connection structure and update rule~\cite{bovier1992rigorous, krotov2016dense, krotov2023new}.

While the connection structure in simplified neuron models is typically assumed to be random or fully connected, spatial networks are also common (and arguably more bio-realistic)~\cite{duchemin2023random,bullmore2009complex, barthelemy2011spatial}. Varying network structure can qualitatively change the dynamics; for example, in the $k$-cap process, spatial network structure leads to stark macroscopic organization in firing patterns~\cite{reid2023k}. In the $N\rightarrow \infty$ limit, spatial networks give rise to local and symmetric continuum dynamics, as in Equation \ref{eq:dynamics}. 

\paragraph{The Method of Moving Planes}The core invariant in our analysis of Equation~\ref{eq:dynamics} follows the \textit{method of moving planes}, a geometric construction in which the reflection of a surface across a family of tangent hyperplanes is used to prove symmetry properties. The method was first used in the original proof of the Alexandrov Soap Bubble Theorem~\cite{alexandrov1962uniqueness}, which shows that, given a bounded connected domain $\Omega\subset \R^n$ with boundary $\Gamma = \partial \Omega$ of class $C^2$, $\Gamma$ is a sphere if and only if it has constant mean curvature. Quantitative variants of the method have been applied to prove partial or near symmetry in a variety of applications~\cite{ciraolo2018sharp,julin2023quantitative}. For a comprehensive overview of the method and its numerous applications, see~\cite{ciraolo2018method}. 

\section{Main Result}
\label{sec:main}
We will prove that the sequence of functions $\psi_t$ become closer to radially symmetric as $t\rightarrow \infty$. The exact structure of the function will depend on the sequence of monotone sharpening functions $\gamma_t$. Formally, we prove the following theorem:

\begin{restatable}{theorem}{thmgeneral}
\label{thm:general}
Let $L_\theta(\psi) = \{x \in \R^d : \psi(x) \ge \theta\}$ denote the superlevel sets of $\psi$. There exists a $p\in \R^d$ such that the following holds. 

For all $\epsilon \in (0, 1)$, there exists a $t^*>0$ such that for all $t > t^*$, and for any $\theta\in (0, 1]$:

\[B(p, \max\{r_\theta - \epsilon, 0\}) \subseteq L_\theta(\psi_t) \subseteq B(p, r_\theta+\epsilon)\]

where $r_\theta$ is the radius of a $d-$ball with volume $|L_\theta(\psi_t)|$

\end{restatable}

For the special case of the fixed-volume $\alpha$-cap process (Definition~\ref{alg:cts}), we show that $A_t$ converges to a ball in $\R^d$ as $t\rightarrow \infty$. 

\begin{figure}
    \centering

\begin{tikzpicture}
%Draw original shape to left
 \begin{scope}
        \filldraw[blue!40, thick] 
            plot[smooth, domain=0:1, samples=100, variable=\x] ({-7.5-5*\x}, {5*sqrt(sqrt(\x) - \x)}) % Upper branch
            -- plot[smooth, domain=1:0, samples=100, variable=\x] ({-7.5-5*\x}, {-5*sqrt(sqrt(\x) - \x)}) % Lower branch (reverse order)
            -- cycle; % Ensures proper closure
    \end{scope}
    \draw[dashed] (-10,-2.7) -- (-10,2.7) node[above] {$H$};
    % Guitar pick shape
    \begin{scope}
        \filldraw[blue!40, thick] 
            plot[smooth, domain=0:1, samples=100, variable=\x] ({-5*\x}, {5*sqrt(sqrt(\x) - \x)}) % Upper branch
            -- plot[smooth, domain=1:0, samples=100, variable=\x] ({-5*\x}, {-5*sqrt(sqrt(\x) - \x)}) % Lower branch (reverse order)
            -- cycle; % Ensures proper closure
    \end{scope}

    % A_t^-
    \begin{scope}
        \filldraw[red!40, thick] 
            plot[smooth, domain=0.5:1, samples=100, variable=\x] ({-5*\x}, {5*sqrt(sqrt(\x) - \x)}) % Upper branch
            -- plot[smooth, domain=1:0.5, samples=100, variable=\x] ({-5*\x}, {-5*sqrt(sqrt(\x) - \x)}) % Lower branch (reverse order)
            --cycle;
 % Ensures proper closure
    \end{scope}
    %A_t^+ \cap A_t^-^H
    \begin{scope}
        \filldraw[green!40, thick] 
            plot[smooth, domain=0.5:1, samples=100, variable=\x] ({-5+5*\x}, {5*sqrt(sqrt(\x) - \x)}) % Upper branch
            -- plot[smooth, domain=1:0.5, samples=100, variable=\x] ({-5 + 5*\x}, {-5*sqrt(sqrt(\x) - \x)}) % Lower branch (reverse order)
            --cycle;
 % Ensures proper closure
    \end{scope}
    \node at (-11.5,0) {$A_t^-$};
    \node at (-8.5,0) {$A_t^+$};
    
    \node at (-4,0) {$A_t^-$};
    \node at (-1,0) {$(A_t^-)^H$};
    \draw (-0.65, 2)-- (-0.5,2.7) node[above] {$A_t^+ \setminus (A_t^-)^H$};
    % Maximal reflecting Hyperplane
    \draw[dashed] (-2.5,-2.7) -- (-2.5,2.7) node[above] {$H$};
\end{tikzpicture}
\caption{(Left) A set $A_t$, with a reflecting hyperplane $H$. (Right) Illustration showing that $H$ is reflecting; $(A_t^-)^H \subseteq A_t^+$.}
\label{fig:Atexample}
\end{figure}
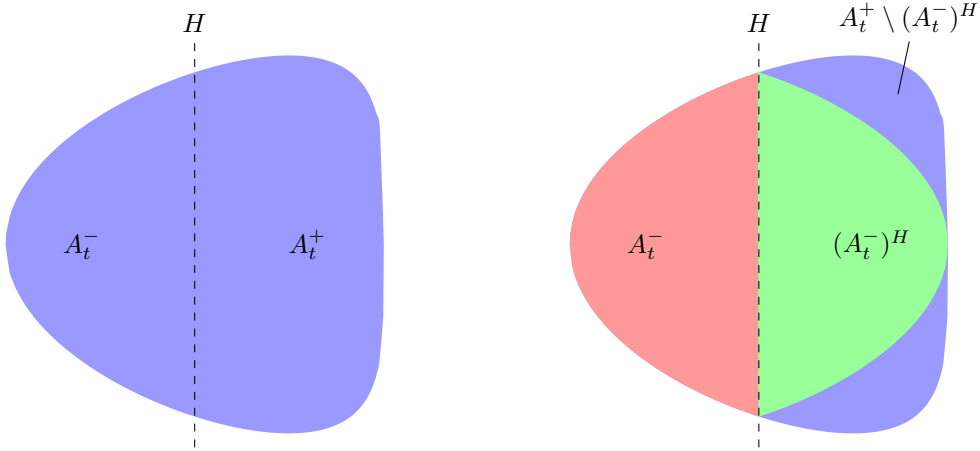

\begin{restatable}{theorem}{thmmain}
\label{thm:main}
    Let $r_*$ be the radius of a ball with volume $|A_0|$. There exists a $p\in \R^d$ such that the following holds. For any $\epsilon>0$, there exists a $t^*$ where for all $t\ge t^*$,
    \begin{align*}
        B(p,r_*-\epsilon)\subseteq A_t\subseteq B(p,r_*+\epsilon)
    \end{align*}
\end{restatable}
The proof of Theorem \ref{thm:general} relies on a key insight, which we refer to as the {\em Reflection Principle}. The operation of the dynamics --- convolution by a symmetric kernel, followed by composition --- preserves reflection symmetry. We can take this idea a step further. Suppose that there exists a hyperplane $H$ such that for all $x\in H^-$, $\psi_t(x) \leq \psi_t(x^H)$. Then, this property remains true inductively for $\psi_{t+1}$. 
\begin{lemma}[Reflection Principle]

\label{lem:reflectionprinciple}
Let $H$ be a hyperplane. Suppose for some $t\ge 1$, for all $x \in H^-$,
\begin{align*}
    \psi_t(x) \le \psi_t(x^H)
\end{align*}
Then, for all $s\ge t$, $\psi_{s}(x) \le \psi_{s}(x^H)$. 
\end{lemma}
\begin{proof}By induction, it suffices to show that the condition $\psi_{t}(x) \le \psi_{t}(x^H)$ implies $\psi_{t+1}(x) \le \psi_{t+1}(x^H)$.
    
    Suppose $\psi_t(x) \le \psi_t(x^H)$. 
    \[f_{t+1}(x) = \int_{\R^d} \psi_t(y)g(x-y)\, dy= \int_{H^-} \psi_t(y)g(x-y)\, dy+ \int_{H^+} \psi_t(y)g(x-y)\,dy \]
    (The integral over the hyperplane $H$ can be ignored as it has measure zero.) Using the reflection notation, $\int_{H^+} \psi_t(y)g(x-y)\, dy = \int_{H^-} \psi_t(y^H)g(x-y^H)\, dy$. Therefore, we can rewrite this:
    \[f_{t+1}(x) = \int_{H^-} \psi_t(y)g(x-y)\, + \psi_t(y^H) g(x - y^H)\, dy\]
    Now, we can examine the difference $f_{t+1}(x^H) -f_{t+1}(x)$:
    \begin{align*}f_{t+1}(x^H)& -f_{t+1}(x) =\\
    &= \int_{H^-} \psi_t(y)g(x^H-y)\, + \psi_t(y^H) g(x^H - y^H) \psi_t(y)g(x-y)- \psi_t(y^H) g(x - y^H)\, dy\\
    &= \int_{H^-} \psi_t(y)(g(x^H-y)- g(x-y)) - \psi_t(y^H)(g(x - y^H)-g(x^H-y^H))\, dy\\
    &=\int_{H^-} (\psi_t(y)- \psi_t(y^H))(g(x - y^H)-g(x^H-y^H))\, dy\\
    \end{align*}
    The last step comes from the fact that $g(x) = \tilde{g}(\|x\|)$ is a function of the $L_2$ norm of $x$; thus, since reflection preserves distance, $g(x-y) = g(x^H - y^H)$.
    
    By the induction hypothesis, $\psi_t(y) - \psi_t(y^H) \le 0$. Further, since $x \in H^-$, $\|x-y^H\| \ge \|x^H - y^H\| \forall y \in H^-$. This implies $g(x - y^H) - g(x^H-y^H) \le 0$. Therefore, the integrand is non-negative, and $f_{t+1}(x^H) \ge f_{t+1}(x)$

    Since $\psi_{t+1}$ is derived by composing $f_{t+1}$ with a monotone function $\gamma_{t+1}$, $\psi_{t+1}(x^H) \ge \psi_{t+1}(x)$ for all $x\in H^-$. Thus, the property holds inductively. 
\end{proof}
\noindent We refer to hyperplanes for which this condition holds as {\em reflecting hyperplanes}. An example of a reflecting hyperplane is illustrated in Figure~\ref{fig:Atexample}. 
\begin{definition}[Reflecting Hyperplane]
A reflecting hyperplane (for $\psi_t$) is a hyperplane $H$ such that $\psi_t(x) \le \psi_t(x^H)$ for all $x \in H^-$. 
\end{definition}

This principle allows us to apply a variant of the \textit{method of moving planes} to argue that $\psi_t$ becomes increasingly radially symmetric with each time step. The argument proceeds as follows. 

First, since $\psi$ is assumed to be bounded, for any fixed direction $u$ there exists a constant $a$ such that $H = \{x : \langle u, x\rangle = a\}$ is a reflecting hyperplane for $\psi$. The argument for this is simple; if $\psi(x) = 0$ for all $x \in H^-$, then $\psi(x) \le \psi(x^H)$ holds trivially. By the boundedness of the support of $\psi$, we can ensure $\psi(x) = 0$ by taking $a$ sufficiently small. Further, it follows that there exists an $a$ such that for any $a' < a$, $H' = \{x : \langle u, x\rangle = a'\}$ is also a reflecting hyperplane for $\psi$. 

Second, if $\psi(x) = \tilde{\psi}(\|x-a\|)$ is a monotone decreasing function of the distance of $x$ to a point $a$, then {\em any} hyperplane $H$ is a reflecting hyperplane for $\psi$ in one of its two orientations.\footnote{Without loss of generality, let $\psi(x) = \tilde{\psi}(\|x\|)$, $u = \mathbf{e}_1$ and $a<0$. Suppose $x\in H^-$. Then $\|x^H\| = \| x- 2(\langle x, u\rangle - a)u\|=\sqrt{(x_1 - 2(x_1-a))^2 + \sum_{i=2}^d x_i^2}$. Since $x\in H^-$, $x_1 < a<0$, and $|x_1 - 2(x_1 - a)| = |-x_1 + 2a|< |x_1|$. This implies $\| x^H \| < \|x\|$, so $\psi(x^H) \ge \psi(x)$.} As we will argue later, an inverse and approximate relationship holds; that is, functions for which `most' hyperplanes are reflecting (under a particular definition of `most') will be close to radially symmetric. 

Hence, by showing that the set of hyperplanes that are non-reflecting is decreasing at each time step, we will argue convergence to radial symmetry. 
\paragraph{Proof outline}The proof of Theorem \ref{thm:general} relies on the reflection principle to establish that the function $\psi_t$ becomes more symmetric in every vector direction with each iteration. The main steps are as follows: 

\begin{itemize}
\item We introduce the concept of a \textit{reflecting slab} in a direction $u$, defined as a set $\{x:a\le \langle x,u\rangle\le b\}$ such that all parallel hyperplanes outside this set are reflecting for $u$ or $-u$. (Definition \ref{def:reflectingslab})

\item We prove that, for every vector $u$, the width $d_t(u)$ of the corresponding reflecting slab strictly decreases with each iteration. (Lemma \ref{lem:reflection_progress})
\item Using Lemma \ref{lem:reflection_progress}, we then establish that $d_t(u) \rightarrow 0$ as $t \rightarrow \infty$.
\item Finally, we argue that the condition $\sup_{u\in\mathbb{S}^{d-1}} d_t(u)\rightarrow 0$ implies convergence in radial symmetry, as stated in Theorem \ref{thm:general}.
\end{itemize}
\subsection{Bounding the Progress in Reflection Symmetry}
\label{sec:progress}
In this section, we will show that $\psi_t$ simultaneously becomes `more symmetric' in every vector direction. The degree of symmetry in a vector direction $u$ will be defined by what we call a \textit{reflecting slab} in the direction of $u$.  

\begin{definition}[Reflecting Slab]
\label{def:reflectingslab}
Let $u\in \mathbb{S}^{d-1}$, $\psi_t:\R^d \rightarrow [0,1]$.

    We call $\{x:a\le \langle x,u\rangle\le b\}$ a reflecting slab for $\psi_t$ (in the direction of $u$) of width $b-a$ if,
    \begin{itemize}
        \item For all $c<a$, $\{\langle x,u\rangle=c\}$ is a reflecting hyperplane.
        \item For all $c>b$, $\{\langle x,-u\rangle=-c\}$ is a reflecting hyperplane.
    \end{itemize}
\end{definition}
The metric which we will track across time steps is the minimum width of a reflecting slab in a particular direction. 
\begin{definition}\label{def:reflecting_slab_width}
    Let $u\in \mathbb{S}^{d-1}$, $\psi_t:\R^d\rightarrow[0,1]$. Define the supremum over reflecting hyperplanes in the direction $u$:
    \[a_t(u) = \sup\{a\in \mathbb{R}:\forall c< a,\; H=\{\langle x, u\rangle = c\} \text{ is a reflecting hyperplane for } \psi_t\}\]

    \noindent Define the smallest width of a reflecting slab in the direction $u$:
    \[d_t(u) = -a_t(-u) - a_t(u)\]
    
    \noindent Finally, let $d_t:=\sup_{u\in\mathbb{S}^{d-1}}d_t(u)$.

    \begin{figure}
    \centering
    \begin{tikzpicture}
        \fill[fill=blue!40] (-2,0) -- (1,1) -- (0,2.25) -- cycle;
        \draw[dashed] (-2,-0.5) -- (-2,3);
        \draw[dashed] (0,-0.5) -- (0,3);
        \fill (0,0) circle (0.1);
        \draw[->] (0,0) -- node[auto] {$u$} (1,0);
        \draw[<->] (-2,2.5) -- node[auto] {$d_t(u)$} (0,2.5);
        \node at (-0.25,1.25) {$A_t$};
        \node at (-2,-0.75) {$a_t(u)$};
        \node at (0,-0.75) {$-a_t(-u)$};
    \end{tikzpicture}
    \caption{Illustration of Definition \ref{def:reflecting_slab_width} in the case where $\psi_t=\mathbbm{1}_{A_t}$.}
    \label{fig:helper1}
    \end{figure}
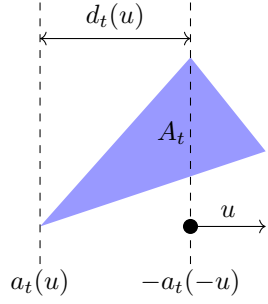
\end{definition}
By the reflection principle and the assumption that the support of $\psi_t$ is bounded, the parameter $a_t(u)$ exists and is finite for all $t$. 

The main lemma of this section, Lemma~\ref{lem:reflection_progress}, will show that $a_t(u)$ is strictly increasing in $t$ for any direction $u\in \mathbb{S}^{d-1}$. This is the most technical lemma of the result, and it involves bounding the gradient of $f_t$ in the direction of $u$ near the hyperplane.
\begin{lemma}
    \label{lem:reflection_progress}
    
    Let $u \in \mathbb{S}^{d-1}$, and let $H=\{ x : \langle x, u \rangle = a\}$ be a reflecting hyperplane. Suppose $H_s=\{ x : \langle x, u \rangle = a+s\}$ is a reflecting hyperplane for $\psi_t$ for all $s\le 0$. 
    
    \noindent Suppose $\Delta = \Delta(H)=\int_{H^-} \psi_t(x^H) - \psi_t(x)\, dx>0$. Fix $R>0$ with the property that $\psi_t(x) = 0$ for all $\|x\| \ge R$.

    \noindent Then, there is some $\epsilon=\epsilon(\Delta,R)>0$ (increasing in $\Delta$) such that for all $s<\epsilon$, the hyperplane
    \begin{align*}
        H_{s}:=\{\langle x,u\rangle=a + s\}
    \end{align*}
    is a reflecting hyperplane for $\psi_{t+1}$.
    
\end{lemma}
\begin{proof}[Proof of Lemma \ref{lem:reflection_progress}] Given that, for all $s\le 0$, $H_s = \{x : \langle x, u\rangle = a+s\}$ is a reflecting hyperplane for $\psi_t$, $H_s$ is also a reflecting hyperplane for $\psi_{t+1}$ by the reflection principle. Thus, the aim of the proof is to demonstrate that $H_s$ is a reflecting hyperplane for $\psi_{t+1}$ whenever $0 < s < \epsilon(\Delta, R)$. 

Recall that $\psi_{t+1} = \gamma_t \circ f_t = \gamma_{t}\circ (\psi_t * g)$. Since $\gamma_t$ is assumed to be monotone increasing, $\psi_{t+1}(x) \ge \psi_{t+1}(y)$ if $f_t(x) \ge f_t(y)$. Hence, $H_s$ is a reflecting hyperplane for $\psi_{t+1}$ if it is reflecting for $f_t$.

To show that $H_s$ is a reflecting hyperplane for $f_t$, we argue that the function $f_t$ is locally monotonic near $H$; that is, there exists an $\eta>0$ such that $f_t(x+s_1u) > f_t(x+s_2u)$ for all  $-\eta\le s_1\le s_2\le \eta$ and $\langle x,u\rangle=a$. Then, we separately prove the following two pieces. 
\begin{claim}
\label{claim:part1} For any $x\in B(0, R)$ such that $a-\eta/4\le \langle x,u\rangle\le a+s$, $f_t(x^{H_s})\ge f_t(x)$.
\end{claim}

\begin{claim}
\label{claim:part2} For any $x\in B(0, R)$ such that $\langle x, u\rangle < a-\eta/4$, $f_t(x^{H_s})\ge f_t(x)$.
\end{claim}

The above two claims combined cover the full space $x\in B(0, R)$; together with a lower bound on $\eta$, this would prove the lemma statement. 

Throughout the proof, we write $f=f_t$, dropping the subscript for clarity. Let $\nabla f(x)$ be the gradient of $f$, and define $\nabla_uf(x) = u \cdot \nabla f(x)$. Define $\alpha = \int_{\R^d} \psi_t(x)\, dx$ (defined under the assumption that the support of $\psi_t$ is bounded). 

Further, we define the following two constants, which exist under the assumption that $g$ is bounded and has continuous first and second derivatives:
\[M = \sup_{x \in [0, 2R]} |\tilde{g}''(x)|\hspace{2.0cm}L = \sup_{x, y \in [0, 2R]}\frac{|\tilde{g}(x) - \tilde{g}(y)|}{|x-y|}\]

\paragraph{Proof of Claim \ref{claim:part1}}First, we will show that $f(x+su)$ is monotone in $s$ in an area near the hyperplane $H$; that is, there exists an $\eta > 0$ such that for $x\in B(0,R)\cap H$ and $-\eta\le s_1\le s_2\le \eta$, $f(x+s_1u)\le f(x+s_2u)$. This will follow if we can show that for any $x\in B(0, R) \cap H$ and $s\in (-\eta, \eta)$, $\nabla_u f(x+su) > 0$. 
    We can write the gradient as follows, applying integration coordinate-wise:
    \begin{align*}
    \nabla f(x+su) &= \int_{\R^d}\psi(y)\nabla g(x-y+su) \, dy\\
    &=\int_{H^+}\psi(y)\nabla g(x-y+su) + \psi(y^H) \nabla g(x-y^H+su)\, dy\\
    \end{align*}
    By Assumption \ref{kernel-assumptions}, $\nabla g(x) = \frac{x}{\|x\|}\tilde{g}'(\|x\|)$ for a monotone decreasing function $\tilde{g}$. Since $x\in H$, for any $y\in \R$, $\|x-y\| = \|x-y^H\|$. Further, since $H$ is orthogonal to $u$, we have that $\langle x-y, u\rangle = -\langle x-y^H, u\rangle$. This gives the following relation to the gradient:
    \[\nabla_u g(x-y^H) = \frac{\langle x-y^H, u\rangle}{\|x-y^H\|}\tilde{g}'(\|x-y^H\|) = \frac{-\langle x-y, u\rangle}{\|x-y\|}\tilde{g}'(\|x-y\|) = -\nabla_u g(x-y)\]
    We will use this to rewrite the gradient bound, as follows. 
\begin{align*}
\nabla_u f(x+su)
&=\int_{H^+}\psi(y)\nabla_u g(x-y+su) + \psi(y^H) \nabla_u g(x-y^H+su)\, dy\\
\\[1pt]
&=\int_{H^+}\psi(y)\nabla_u g(x-y) + \psi(y^H) \nabla_u g(x-y^H)\\
&\mathrel{\phantom{=}} + \psi(y^H)\bigl(\nabla_u g(x-y^H+su)-\nabla_u g(x-y^H)\bigr)\\
&\mathrel{\phantom{=}} + \psi(y)\bigl(\nabla_u g(x-y+su)-\nabla_u g(x-y)\bigr)\, dy\\
\\[1pt]
&= \int_{H^+}(\psi(y)-\psi(y^H))\nabla_u g(x-y)\\
&\mathrel{\phantom{=}} + \psi(y^H)\bigl(\nabla_u g(x-y^H+su)-\nabla_u g(x-y^H)\bigr)\\
&\mathrel{\phantom{=}} + \psi(y)\bigl(\nabla_u g(x-y+su)-\nabla_u g(x-y)\bigr)\, dy
\end{align*}
    Under the assumption that $|\tilde{g}''(x-y)|\leq M$ whenever $x, y \in B(0, R)$, we can bound the sum of the last two terms above:
    \begin{align*}\int_{H^+}& \psi(y^H)\bigl(\nabla_u g(x-y^H+su)-\nabla_u g(x-y^H)\bigr) + \psi(y)\bigl(\nabla_u g(x-y+su)-\nabla_u g(x-y)\bigr)\, dy\\
    &\geq\int_{H^+}-\psi(y^H)|s|M - \psi(y)|s|M\, dy \\
    &\geq -\alpha \eta M
    \end{align*}
    
    If the first term of the integral is bounded above 0 by a constant, then we can choose $\eta$ sufficiently small such that the gradient is positive. However, this term requires a more technical argument to bound. Writing out the gradient:
    
    \begin{equation}   
    \int_{H^+}(\psi(y)-\psi(y^H))\nabla_u g(x-y)= \int_{H^+}(\psi(y)-\psi(y^H))\frac{\langle x-y, u\rangle}{\|x-y\|}\tilde{g}'(\|x-y\|)
    \label{eq:integralbound}
    \end{equation}

    Note that for $y \in H^+$, $\langle y, u\rangle > \langle x, u\rangle=a_t(u)$ and $\psi(y) \geq \psi(y^H)$. Also, $\tilde{g}'(c)<0$ for all $c>0$. This implies that the integral is non-negative. However, the integrand can be arbitrarily near zero if $y$ is near the hyperplane. Similarly, $\tilde{g}'(\|x-y\|)/\|x-y\|$ may approach zero as $y$ approaches $x$. Therefore, we will employ a geometric argument to bound the integral away from zero,. Define \[S_{\delta, r} = \{y: a\le \langle u, y\rangle < a+\delta, y \in B(0, R)\} \cup \{y : \lVert x-y\rVert \leq r\}\] 
    This set represents the points $y$ where the integrand of Equation \ref{eq:integralbound} is near zero, and it is illustrated in Figure~\ref{fig:geometric-arg}. The volume of $S_{\delta, r}$ can be bounded, letting $V_d$ be the volume of the unit $d-$ball. 
    \begin{align*}
    \vol(S_{\delta, r}) &\le \vol(B(0, R) \cap \{y: a\le \langle u, y\rangle < a+\delta\}) + \vol(\{y : \lVert x-y\rVert \leq r\})\\
    &\le \int_{a}^{a+\delta}\vol(B(0, R) \cap \{\langle u, y\rangle = s\})\,d s  + r^d V_d \\
    &\le \int_a^{a+\delta} V_{d-1}R^{d-1}\, ds + r^d V_d\\
    &= \delta V_{d-1}R^{d-1} + r^d V_d
    \end{align*} 
    Setting $\delta = \frac{\Delta}{4V_{d-1} R^{d-1}}$ and $r = \left(\frac{\Delta}{4V_d}\right)^{1/d}$, the volume of this region is at most $\frac{\Delta}{2}$. 

 Define the following constant\footnote{As a sanity check, note that $\Delta \leq \text{Vol}(B(0, R)) = V_d R^d$, so the interval is nonempty.}:
 \[\lambda_{\Delta, R} \coloneqq \inf\left\{\dfrac{|\tilde{g}'(c)|}{c}: c\in \left(\left(\frac{\Delta}{4V_d}\right)^{1/d}, 2R\right]\right\}\]
 Since it is assumed that $\tilde{g}'$ is continuous and $|\tilde{g}'(x)|>0$ when $x>0$, this constant is bounded away from zero. This leads to the following inequality, letting $S_x = \{y\in H^+: \langle u,y\rangle \geq a+\delta, \lVert x-y \rVert \geq r\}$:
    
    \begin{align*}\nabla_u f(x+su) &\geq \int_{S_x}(\psi(y)-\psi(y^H))\frac{\langle x-y, u\rangle}{\lVert x-y\rVert}\tilde{g}'(\lVert x-y\rVert)\, dy - \alpha\eta M\\
    &\geq \int_{S_x}(\psi(y)-\psi(y^H)) \lambda_{\Delta, R}\delta\, dy - \alpha \eta M\\
    &\geq \frac{\Delta\delta }{2}\lambda_{\Delta, R}-\alpha\eta M 
    \end{align*}
    since $0<\psi(y) - \psi(y^H) < 1$. 

    Therefore, by picking $\eta = \frac{\Delta \lambda_{\Delta, R} \delta}{2\alpha M}$, $\nabla_u f(x+su) > 0$ for any $s\in (-\eta, \eta)$; and, for  $-\eta\le s_1\le s_2\le \eta$, $f(x+s_1u)\le f(x+s_2u)$.
\begin{figure}
\centering
\begin{tikzpicture}
    
    % Define parameters
    \def\a{-2.5} % x = a
    \def\del{1} % delta width
    \def\RR{3} % Large ball radius
    \def\r{1.5} % Small ball radius
    \def\x{-2.5} % x-coordinate of x
    \def\y{0} % y-coordinate of x

 \begin{scope}
        \filldraw[blue!40, thick] 
            plot[smooth, domain=0:1, samples=100, variable=\x] ({-5*\x}, {5*sqrt(sqrt(\x) - \x)}) % Upper branch
            -- plot[smooth, domain=1:0, samples=100, variable=\x] ({-5*\x}, {-5*sqrt(sqrt(\x) - \x)}) % Lower branch (reverse order)
            -- cycle; % Ensures proper closure
    \end{scope}

    % A_t^-
    \begin{scope}
        \filldraw[gray!40, thick] 
            plot[smooth, domain=0.5:1, samples=100, variable=\x] ({-5*\x}, {5*sqrt(sqrt(\x) - \x)}) % Upper branch
            -- plot[smooth, domain=1:0.5, samples=100, variable=\x] ({-5*\x}, {-5*sqrt(sqrt(\x) - \x)}) % Lower branch (reverse order)
            --cycle;
 % Ensures proper closure
    \end{scope}
    %A_t^+ \cap A_t^-^H
    \begin{scope}
        \filldraw[gray!40, thick] 
            plot[smooth, domain=0.5:1, samples=100, variable=\x] ({-5+5*\x}, {5*sqrt(sqrt(\x) - \x)}) % Upper branch
            -- plot[smooth, domain=1:0.5, samples=100, variable=\x] ({-5 + 5*\x}, {-5*sqrt(sqrt(\x) - \x)}) % Lower branch (reverse order)
            --cycle;
 % Ensures proper closure
    \end{scope}
    % Shade right half-plane H^+ inside B(0, R), excluding the blue strip and red ball

    % Remove (cut out) the strip and the small ball by overlaying white
    \begin{scope}
        \clip (\a, -\RR) rectangle (\RR, \RR); % Restrict clipping to H^+
        \fill[red!30] (\a,-\RR) rectangle (\a+\del,\RR); % Cut out the blue strip
        \fill[red!30,] (\x,\y) circle (\r); % Cut out the small red ball
    \end{scope}
    % Draw large ball B(0, R)
    \draw[dashed] (\x,\y) circle (\r);
    % Draw vertical line H
    \draw[dashed] (\a,-\RR) -- (\a,\RR) node[above] {$H$};

    % Mark points
    \fill (\x,\y) circle (2pt) node[below left] {$x$};
    
    \draw (\x, \y) -- (-1.2, -0.8) node[pos=0.5, above]{$r$};
    \draw (\x, 2) -- (\x + \del, 2) node[pos=0.5, above]{$\delta$};

    %\node at (\x + 2*\del, \y+1.5){$S_x$};
\end{tikzpicture}

\caption{Illustration of the geometric argument in bounding $\nabla_u f(x)$, using the example $\psi_t(x) = \mathbbm{1}\{x\in A_t\}$ from Figure \ref{fig:Atexample}. The volume of the red region is at most $\Delta/2$. The blue region indicates the points $y$ such that $\psi_t(y) > \psi_t(y^H)$.}
\label{fig:geometric-arg}
\end{figure}
Let $s\in[0,\eta/4]$. Then, for $x\in B(0,R)$ with $a-\eta/4\le \langle x,u\rangle\le a+s$, $f(x^{H_s})\ge f(x)$. This implies Claim \ref{claim:part1}.

\paragraph{Proof of Claim \ref{claim:part2}}Now, we show that $H_s$ is reflecting for all of $B(0, R)$ by arguing that for $x : \langle x, u\rangle<a-\eta/4$, $f(x^{H_s}) \ge f(x)$.

    Let
    \begin{align*}
        h(s)=\inf\{f(x^{H_s})-f(x) \mid x\in B(0,R)\cap \{\langle x,u\rangle \le a-\eta/4\}\}
    \end{align*}
    Note that if $h(s)\ge 0$, then for $x\in B(0,R)$ with $\langle x,u\rangle\le a-\eta/4$, $f(x^{H_s})\ge f(x)$. Combined with Claim \ref{claim:part1}, this would imply that $H_s$ is a reflecting hyperplane for $\psi_t$. 
    
    By the assumption that $g$ is $L$-Lipschitz (Assumption \ref{kernel-assumptions}), $f$ is $\alpha L$-Lipschitz. Since $\|x^{H_s}-x^{H_t}\|=2|s-t|$, we have $|f(x^{H_s})-f(x^{H_t})|\le 2\alpha L|s-t|$. Thus, $f(x^{H_s})-f(x)$ is $2\alpha L$-Lipschitz with respect to $s$ for any fixed $x$. Since $h(s)$ is the infimum of such functions, $h$ is $2\alpha L$-Lipschitz as well.

    Since $f(x^H) \geq f(x)$ by the reflection principle, $h(0)\ge 0$. We can in fact lower bound it by a positive constant depending only on $R$ and $\Delta$.
        
    The rotation symmetry of $g$ implies that $g(x-y) = g(x^H-y^H)$, allowing us to write $f(x^H) - f(x)$ as follows:
    \begin{align*}
f(x^H) - f(x) &= \int_{\R^d}\psi(y)g(x^H-y) -\psi(y) g(x-y)\, dy\\
&= \int_{H^+}\psi(y)g(x^H-y) -\psi(y) g(x-y)\\
&\mathrel{\phantom{=\int_{H^+}}} +\,\psi(y^H)g(x^H-y^H) -\psi(y^H) g(x-y^H)\, dy\\
&= \int_{H^+}\psi(y)g(x^H-y) -\psi(y) g(x-y)\\
&\mathrel{\phantom{=\int_{H^+}}} +\,\psi(y^H)g(x-y) -\psi(y^H) g(x^H-y)\, dy\\
&=\int_{H^+}\bigl(\psi(y)-\psi(y^H)\bigr)\bigl(g(x^H-y)-g(x-y)\bigr)\, dy\\
    \end{align*}
    
    By the assumption that $H$ is reflecting, this integral is nonnegative. 
    
    To give a nontrivial lower bound, a similar geometric argument is again necessary. If $y$ is too close to the hyperplane, then the differences in distance $\lVert x-y\rVert$ and $\lVert x^H - y\rVert$ can be arbitrarily small. So, as above, let $S_{x^H} = \{y\in H^+: \langle u,y\rangle \geq a+\delta, \lVert x^H-y \rVert \geq  \left(\frac{\Delta}{4V_d}\right)^{1/d}\}$, using the same constant $\delta= \frac{\Delta}{4V_{d-1} R^{d-1}}$. Bounding the difference with the integral over this set: 
    \[f(x^H) - f(x) \geq \int_{S_{x^H}} (\psi(y)-\psi(y^H))(g(x^H - y) - g(x-y))\, dy\]

The integrand can be bounded by a function of $\eta, \delta$, and $\lambda_{\Delta, R}$. For any $x$ with $\langle x, u\rangle < a-\eta/4$ and $y\in S_{x^H}$, let $\rho=\|x^H - y\|$ and $\sigma=\|x-y\|$.
\begin{align*}
\tilde g(\rho)-\tilde g(\sigma) &\ge 
   \tilde g\!\left(\frac{\rho+\sigma}{2}\right)-\tilde g(\sigma) \\[4pt]
&\ge 
   \left(\sigma-\frac{\rho+\sigma}{2}\right)
   \min_{c\in[(\rho+\sigma)/2,\ \sigma]}
        \tilde g'(c) \\[4pt]
&\ge 
   \frac{(\sigma-\rho)(\rho+\sigma)}{4}\,
   \min_{c\in[(\rho+\sigma)/2,\ \sigma]}
        \frac{\tilde g'(c)}{c} \\[4pt]
&\ge 
   \frac{(\sigma-\rho)(\rho+\sigma)}{4}\,
   \min_{c\in[(\Delta/4V_{d})^{1/d},\ 2R]}
        \frac{\tilde g'(c)}{c} \\[4pt]
&=
   \frac{\lambda_{\Delta,R}}{4}\,
      (\sigma^{2}-\rho^{2}),
\end{align*}
The square difference  $\sigma^2 - \rho
^2 = \|x-y\|^2 - \|x^H - y\|^2$, for $y\in S_{x^H}$, is bounded as follows:

\begin{align*}
\|x- y\|^2-\|x^H-y\|^2 &=\|x-x^H\|^2 + 2\langle x-x^H, x^H - y\rangle\\
&= \langle x-x^H, x+x^H - 2y\rangle\\
&=2\langle x^H-x, y- 1/2(x+x^H)\rangle\\
&=2\langle 2(a-\langle x, u\rangle)u, y- 1/2(x+x^H)\rangle = 4(a-\langle x, u\rangle)\langle u, y\rangle\\
&\ge \eta \delta
\end{align*}
This leads to the following bound on $h(0)$. For any $x\in B(0, R)$ with $\langle x, u \rangle \le a - \eta/4$

\[f(x^H) - f(x) \geq \int_{S_{x^H}} (\psi(y) - \psi(y^H))g(y-x^H) - g(y-x)\, dy\geq \dfrac{\Delta}{2}\cdot \dfrac{\lambda_{\Delta, R}}{4}\cdot \eta\delta= \frac{\lambda_{\Delta, R}\Delta \eta \delta}{8}\]
This implies $h(0) \ge \frac{\lambda_{\Delta, R}\Delta \eta \delta}{8}$
Finally, we can set the value of $\epsilon(\Delta, R)$ as follows:
    \[\epsilon=\min\left(\dfrac{\eta}{4},\dfrac{h(0)}{2L}\right)=\min\left(\dfrac{\eta}{4}, \dfrac{\lambda_{\Delta, R}\Delta \eta \delta}{8}\right)\: \text{where }\eta=\frac{\Delta \lambda_{\Delta, R} \delta}{2\alpha M}\]
Then, for $0\le s\le \epsilon$, $h(s)\ge h(0)-2L\epsilon\ge 0$. So, for any $x\in B(0,R)$ with $\langle x,u\rangle\le a-\eta/4$, $f(x^{H_{s}})\ge f(x)$. By the first part of the proof (Claim \ref{claim:part1}), for $x\in B(0,R)$ with $a-\eta/4<\langle x,u\rangle \le a+s$, $f(x^{H_s})\ge f(x)$.
    
    Thus, for all $x\in B(0, R)\cap H_s^-$, $f(x^{H_s})\ge f(x)$. By the assumption that the support of $\psi_{t+1}$ is bounded within $B(0, R)$, $\psi_{t+1}(x^{H_s}) \ge \psi_{t+1}(x)$. This is true for all $0 \le s \le \epsilon$, so $a_{t+1}(u) \ge a_t(u) + \epsilon$. 
    
\end{proof}

Lemma \ref{lem:reflection_progress} demonstrates that $d_{t+1}(u) \le d_t(u) - \epsilon_t$, where $\epsilon_t$ is a function of $\Delta_t$, defined as
\[\Delta_t = \lim_{a \rightarrow a_t(u)} \int_{H^-} \psi_t(y^H) - \psi_t(y)\, dy\] 
However, this alone does not imply convergence to radially symmetric function, since a priori $d_t(u)$ may approach a limit other than 0. Next, Lemma \ref{lem:zero_dt} will establish a lower bound on $\Delta_t$ given $d_t$, and use this to argue that $d_t\rightarrow 0$. 

\begin{lemma}
\label{lem:zero_dt}
In fact, $\lim_{t\to\infty}d_t=0$.
\end{lemma}

Toward the proof of this Lemma, we first show the following bound on the total mass of $\psi_t$ contained within a reflecting slab. 
\begin{lemma}
\label{lem:intersectionvolume}
Let $\phi:\R^d\rightarrow [0,1]$ be an integrable function with $\{x : \phi(x) > 0\}\subseteq[0, 2R]^d$, and let $0 \le a < b \le 2R$ with $l=b-a$. Also, let $u\in \mathbb{R}^d$ be any unit vector.

Suppose $H_a = \{\langle x, u\rangle = a\}$ and $H_b = \{\langle x, -u\rangle = -b\}$ are both reflecting hyperplanes for $\phi$. Then, 
\[\int_{H_a^+ \cap H_b^+}\phi(x)\, dx\ge \frac{\int_{[0, 2R]^d}\phi(x)\, dx}{\lceil \frac{2R}{l}\rceil + 1}\]

\end{lemma}
\begin{proof}
Assume without loss of generality that $\int_{[0, 2R]^d}\phi(x)\, dx=1$.

Divide $[0, 2R]^d$ into slabs; let $I_k = \{x\in [0, 2R]^d : a-kl \le \langle x, u \rangle < a-(k-1)l\}$. Note that $\{I_k : k = -\lfloor \frac{2R-a}{l}\rfloor, \dots, -1, 0, 1, \dots, \lceil \frac{a}{l}\rceil\}$ covers all of $[0, 2R]^d$. 

Let $I_m$ be the slab maximizing $\int_{I_k}\phi(x)\, dx$. By the pigeonhole principle, $|S \cap I_m| \ge \frac{1}{\lceil \frac{2R}{l}\rceil + 1}$. 

Note that $I_0$ is $H_a^+ \cap H_b^+$, the slab of interest. We will show by induction that $\int_{I_0}\phi(x)\, dx \ge \int_{I_m} \phi(x)\, dx$.

First, suppose $m>0$. By the assumption that $H_a$ is a reflecting hyperplane, for any $x \in I_m \subseteq H_a^-$, $\phi(x^{H_a})\ge \phi(x)$. By definition, $x^{H_a} = x - 2(\langle x, u \rangle-a) u$. So, 
\[\langle x^{H_a}, u \rangle = \langle x, u\rangle - 2 (\langle x, u\rangle - a)\langle u, u\rangle = 2a-\langle x, u\rangle \]
Since we know  $a-ml \le \langle x, u \rangle < a-(m-1)l$, 
\[a + (m-1)l \le \langle x^{H_a}, u \rangle < a+ml\]
Thus, $x^{H_a} \in I_{-(m-1)}$. This implies that $\int_{I_{-(m-1)}} \phi(x)\, dx\ge \int_{I_m}\phi(x)\, dx$.

Similarly, suppose $m<0$. By the assumption that $H_b$ is a reflecting hyperplane, for any $x \in I_m \subseteq H_b^-$, $\phi(x^{H_b})\ge \phi(x)$. Again, since $x^{H_b} = x - 2(\langle x, u \rangle-b) u$, $\langle x^{H_b}, u \rangle = 2b-\langle x, u\rangle$. Since we know  $a-ml \le \langle x, u \rangle < a-(m-1)l$, 
\[2b-(a-(m-1)l) \le \langle x^{H_b}, u\rangle \le 2b - (a-ml)\]
Applying $l=b-a$, 
\[a + (m+1)l \le \langle x^{H_b}, u\rangle \le a+(m+2)l\]
Thus, $x^{H_b} \in I_{-(m+1)}$. This implies that $\int_{I_{-(m+1)}}\phi(x)\, dx \ge \int_{I_m}\phi(x)\, dx$.

In both cases, the absolute value of the index decreases by 1. By induction, we can say $\int_{I_0}\phi(x)\, dx \ge\int_{I_m}\phi(x)\,dx$. This proves the statement of the lemma. 
\end{proof}

\begin{lemma}
There is a function $\epsilon_2(d,R)>0$ (increasing in $d$) such that for any $u\in\mathbb{S}^{d-1}$,
\begin{align*}
    d_{t+1}(u)\le d_t(u)-\epsilon_2(d_t,R)
\end{align*}
\label{lem:decrease_dt}
\end{lemma}
\begin{proof}
    Let $\alpha = \int_{\R^d}\psi_t(x)\, dx$

    Define the hyperplanes:
    \begin{align*}
        H_a&=\{\langle x,u\rangle=a\}&H_b&=\{\langle x,-u\rangle=-b\}
    \end{align*}
    to be any two reflecting hyperplanes with the property that
    \begin{align*}
        &\{\langle x,u\rangle=c\} \text{ is a reflecting hyperplane for all $c\le a$}\\
        &\{\langle x, -u\rangle =-c\} \text{ is a reflecting hyperplane for all $c\le b$}
    \end{align*}
    By Lemma \ref{lem:intersectionvolume}, we see that $H_a^+\cap H_b^+$ contains at least a $\delta = \dfrac{1}{\lceil \dfrac{2R}{b-a}\rceil+1}$ fraction of the measure of $\psi_t$.
    
    Then, one of $H_b^-$ and $H_a^-$ must have a fraction of the measure of $\psi_t$ at most $\dfrac{1}{2}-\dfrac{\delta}{2}$; WLOG suppose $\dfrac1\alpha\int_{H_a^-}\psi_t(x)\, dx \le \dfrac{1}{2} - \dfrac{\delta}{2}$. Here, $\Delta =\int_{H_a^-} \psi_t(x^H) - \psi_t(x)\, dx\ge \alpha\delta = \dfrac{\alpha}{\lceil \frac{2R}{b-a}\rceil+1}$.
    
    By Lemma \ref{lem:reflection_progress}, taking $\hat{\epsilon} = \epsilon\left(\dfrac{\alpha}{\lceil \frac{2R}{b-a}\rceil+1}, R\right) > 0$, $H_s = \{\langle x, u\rangle  = a+s\}$ is a reflecting hyperplane for $\psi_{t+1}$ for all $s\in [0, \hat{\epsilon})$. Since $H_b$ is still a reflecting hyperplane for $\psi_{t+1}$ by the reflection principle (Lemma \ref{lem:reflectionprinciple}), there exists a reflecting slab at $t+1$ defined by $H_b$ and $H_{a+\hat{\epsilon}}$. 

    By definition, we can take $a$ and $b$ arbitrarily close to $a_t(u)$ and $b_t(u)$; this implies that, e.g., $d_{t+1}(u) = b_{t+1}(u) - a_{t+1}(u) < b_t(u) - (a_t(u) + \hat{\epsilon}/2) = d_t(u) + \hat{\epsilon}/2$. 

    By the definition of $d_t$, there exists a reflecting slab of width $d_t$ in every direction $u$. Therefore, we can set 
    \[\epsilon_2 = \epsilon\left(\dfrac{\alpha}{\lceil \frac{2R}{d_t}\rceil+1}, R\right) > 0\]
    
    which implies that $d_{t+1} \le b - (a+\epsilon_2) = d_t - \epsilon_2/2$. Further, $\epsilon_2$ is increasing in $d_t$. 
\end{proof}

Finally, we can prove the statement on the limit of $d_t$.
\begin{proof}[Proof of Lemma \ref{lem:zero_dt}]
    For every $u\in\mathbb{S}^{d-1}$, $\psi_t$ has a reflecting slab of width $d_t$ in direction $u$.
    By the previous lemma, for every $u\in\mathbb{S}^{d-1}$, $\psi_{t+1}$ has a reflecting slab of width $d_t-\epsilon_2(d_t,R)$ in direction $u$.
    Hence, $d_{t+1}\le d_t-\epsilon_2(d_t,R)$

    By monotone convergence, $\lim_{t\to\infty}d_t=d_*\ge 0$ for some $d_*$. Suppose $d_*>0$. Then, $d_{t+1}-d_t\le -\epsilon_2(d_t,R)\le -\epsilon_2(d_*,R)$. However, since $\epsilon_2(d^*, R)$ is bounded away from zero, this is a contradiction. Hence, the limit must be zero. 
\end{proof}

\subsection{Convergence to Radial Symmetry}
\label{sec:convergence}
In this section, we will prove that as $t\rightarrow \infty$ (and therefore $d_t \rightarrow 0$), the functions $\psi_t$ become arbitrarily close to radially symmetric about a common point $p$. First, we will need the following result on the intersection of reflecting slabs:
\begin{lemma}
\label{lem:slab_intersection_unique}
    The reflecting slabs of $\psi_t$ for $t\ge 0$ have a unique common intersection.
\end{lemma}
The proof of this result is deferred to the appendix. Next, we prove a geometric relationship between the structure of $\psi_t$ and the width of the intersection over reflecting slabs (Lemma \ref{lem:helper0}). Lastly, we will use this lemma to prove the main theorem. 
\begin{lemma}
\label{lem:helper0}
Let $p\in \R^d$, $\epsilon>0$, and $\psi:\R^d\rightarrow[0,1]$ have the property that all hyperplanes $H$ for which $B(p,\epsilon)\subseteq H^+$ are reflecting hyperplanes.

Then, for any $x, y\in \R^d$ satisfying $\|x-p\|-\|y-p\|>2\epsilon$, $\psi(x) \le \psi(y)$.
\end{lemma}
\begin{proof}
Let
\begin{align*}
    H_{x\to y}
    =\left\{v\in\R^d:\left\langle v, \frac{y-x}{\|y-x\|}\right\rangle=\left\langle \frac{x+y}{2}, \frac{y-x}{\|y-x\|}\right\rangle\right\}
\end{align*}
be the hyperplane such that reflecting $x$ across it yields $y$.

We will show that $B(p,\epsilon)\subseteq H_{x\to y}^+$, which would imply that $H_{x\to y}$ is a reflecting hyperplane, and hence $\psi(x)\le \psi(y)$.

    It suffices to show that
    \begin{align*}
        \left\langle \frac{x+y}{2},\frac{y-x}{\|y-x\|}\right\rangle\le \left\langle p,\frac{y-x}{\|y-x\|}\right\rangle-\epsilon
    \end{align*}

    This follows from the assumption that $\|x-p\|-\|y-p\|\ge 2\epsilon$:
    \begin{align*}
        \left\langle p-\dfrac{x+y}{2},\frac{y-x}{\|y-x\|}\right\rangle&=\left\langle\dfrac{(x-p)+(y-p)}{2},\frac{(x-p)-(y-p)}{\|x-y\|}\right\rangle\\
        &=\frac{\|x-p\|^2-\|y-p\|^2}{2\|x-y\|}\\
        &\ge \frac{\|x-p\|^2-\|y-p\|^2}{2(\|x-p\|+\|y-p\|)}\\
        &=\frac{\|x-p\|-\|y-p\|}{2}\\
        &\ge \epsilon
    \end{align*}
    
    \begin{figure}
    \centering
    \begin{tikzpicture}
        \draw[dashed,fill=gray] (0,0) circle (1);
        \draw[dashed] (0,0) circle (2);
        \draw (0,0) circle (4);
        \draw[dashed] (0,0) -- (0,4);
        \node at (-0.2,3) {$2\epsilon$};
        \node at (-0.2,0.5) {$\epsilon$};
        \node[circle,fill,label=left:$p$] (P) at (0,0) {};
        \node[circle,fill,label=above:$x$] (X) at (0,4) {};
        \node[circle,fill,label=below:$y$] (Y) at (1.2,-1.6) {};
        \draw[->] (X) -- (Y);
    \coordinate (M) at ($(X)!0.5!(Y)$);
    \draw (M) node[below right] {$H_{x\to y}$};
    % H_{x\to y} is perpendicular bisector of xy
    \draw[thick] ($(M)!2!90:(X)$) -- ($(M)!2!-90:(X)$);
    \end{tikzpicture}
    \caption{Illustration of Lemma \ref{lem:helper0}. Note that $B(p,\epsilon)\subseteq H_{x\to y}^+$.}
    \label{fig:helper0}
    \end{figure}

\end{proof}

Finally, we have the tools we need to prove the main convergence theorem, restated below for clarity. 

\thmgeneral*
\begin{proof}
Let $p$ be the common intersection of all reflecting slabs for $\psi_t$ for $t \ge 0$ (per Lemma \ref{lem:slab_intersection_unique}).

Fix $\epsilon > 0$. Per Lemma \ref{lem:zero_dt}, there exists a time $t^*>0$ such that for all $t \ge t^*$, $d_{t} \le \frac{1}{3}\epsilon$; that is, the reflecting slab in every direction $u$ has thickness at most $\frac{1}{3}\epsilon$. From Definition \ref{def:reflectingslab}, it follows that any hyperplane $H$ satisfying $B(p,\frac{1}{3}\epsilon)\subseteq H^+$ is a reflecting hyperplane. Thus, we can apply Lemma \ref{lem:helper0}; for any $x, y$ satisfying $\|x-p\|-\|y-p\|>\frac{2}{3}\epsilon$, $\psi_t(x) \le \psi_t(y)$.

Fix $\theta\in (0, 1]$. Let $r_\theta$ be the radius of a ball with volume $L_\theta(\psi_t)$. Since the support of $\psi_t$ is assumed to be bounded, this is well defined. 

We now show the left inclusion, $B(p, \max\{r_\theta - \epsilon, 0\}) \subseteq L_\theta(\psi_t)$.

If $r_\theta<\epsilon$, $\emptyset\subseteq L_\theta(\psi_t)$ holds trivially, so we suppose otherwise.

Since $|L_\theta(\psi_t)|>|B(p,r_\theta-\frac{\epsilon}{3})|$, there exists an $x\in L_\theta(\psi_t)$ for which $\|x-p\|\ge r_\theta-\frac{\epsilon}{3}$. Thus, by Lemma \ref{lem:helper0}, for all $y\in B(p,r_\theta-\epsilon)$, $\psi_t(x) \le \psi_t(y)$. Since $x \in L_\theta(\psi_t)$, we have $y\in L_\theta(\psi_t)$.

We now show the right inclusion, $L_\theta(\psi_t) \subseteq B(p, r_\theta+\epsilon)$. 

Since $|L_\theta(\psi_t)|<|B(p,r_\theta+\frac{\epsilon}{3})|$, there exists $y\notin L_\theta(\psi_t)$ for which $\|y-p\|<r_\theta+\frac{\epsilon}{3}$. Thus, by Lemma \ref{lem:helper0}, for all $x\notin B(p,r_\theta+\epsilon)$,  $\psi_t(x) \le \psi_t(y)$. This implies $x\notin L_\theta(\psi_t)$.

Hence, for all $\theta \in (0, 1]$,
\[B(p, \max(r_\theta - \epsilon, 0)) \subseteq L_\theta(\psi_t) \subseteq B(p, r_\theta+\epsilon)\]
as claimed.
\end{proof}

\section{Volume-Preserving Threshold}
\label{sec:volume-preserving}
In this section, we highlight the special case where $\gamma_t(x) = \mathbbm{1}\{x \ge C_t\}$ for a threshold $C_t\in \R^+$, as given in Definition~\ref{alg:cts}. Let $A_t$ denote the set $\{x : f_t(x) \ge C_t\}$, so $\psi_t = \mathbbm{1}_{A_t}$. We focus on the case where $C_t$ is chosen to preserve the volume of $A_t$ across time steps. This process is equivalent to the MBO scheme for volume-preserving motion by mean curvature~\cite{merriman1992diffusion, ruuth2003simple, laux2017convergence}. We also resolve an open question in this line of work by demonstrating that the threshold $C_t$ is well defined under the condition that $g$ is analytic~\cite{ruuth2003simple, laux2017convergence}. 

It is not immediately obvious whether such a volume preserving threshold exists for a given kernel $g$ and set $A_t$. In Definition~\ref{alg:cts}, we define the threshold $C_t$ as the infimum over all thresholds that yield a set with volume less than the target. In the next lemma, we will show that under the additional assumption that $g$ is an analytic function, this construction ensures that the volume of $A_t$ is preserved. 
\begin{lemma}
\label{lem:analytic}
In addition to Assumption \ref{kernel-assumptions}, suppose that $g$ is analytic. Then, for any compact set $A$ with $|A| = \alpha$, and $f(x) = \int_{A} g(x, y)\, dy$, $\exists c$ with $|\{x : f(x) \ge c\}| = \alpha$.
\end{lemma}
\begin{proof}
Define $h:(0, \infty)\rightarrow\mathbb{R}$ where $h(c) = \text{Vol}(\{x : f(x) \geq c\})$. Since $\{x: f(x) \geq c_1 \} \subseteq \{x: f(x) \geq c_2\}$ whenever $c_1 \geq c_2$, $h$ is monotonically decreasing. 

It will suffice to show that $h$ is continuous on an interval $[C_0,  C_1]$ where $h(C_0) > \alpha$ and $0 \le h(C_1) < \alpha$. The existence of $C_1$ follows from the boundedness of $g$; let $C_1 = \alpha \max_{r} \tilde{g}(r) = C_1$. Then, $f(x) = \int_A \tilde{g}(\|x-y\|)\, dy < C_1$, and $h(C_1) = 0< \alpha$. The existence of $C_0$ also follows from Assumption \ref{kernel-assumptions} and the compactness of $A$. Suppose $A \subset B(0, r)$; then, for all $x \in B(0, 2r)$, $f(x) \ge \alpha \tilde{g}(2r)$. Since $|B(0, 2r)| > |A| = \alpha$, $h(\alpha \tilde{g}(2r)) >\alpha$. 

Suppose for contradiction that $h$ is not continuous on this interval. Then, $\exists c\in [C_0, \infty)$ such that $\lim_{x\rightarrow c^-}h(x) > \lim_{x\rightarrow c^+}h(x)$.
By the monotone convergence theorem,
\begin{align*}
    \lim_{x\to c^-}h(x)&=\text{Vol}(\{x : f(x) \geq c\})=\text{Vol}(\{x : f(x) >c\})+\text{Vol}(\{x : f(x) = c\})\\
    \lim_{x\to c^+}h(x)&=\text{Vol}(\{x : f(x) >c\}
\end{align*}
If $h$ has a discontinuity at $c$, then $\text{Vol}(\{x : f(x) = c\})>0$. We will show that this is a contradiction by demonstrating that $f$ is analytic.

\begin{claim}
 $f(x) = \int_A g(x, y)\, dy$ is analytic on $\R^d$. 
\end{claim}
\begin{proof}
Suppose $A\subset B(0, r)$ by the assumption of compactness. We will prove that $f$ is analytic on $B(0, kr)$ for any $k$. 

Since $g$ is analytic on $\R^d$, there exists an open neighborhood $V\subset \mathbb{C}^d$ such that $\R^d\subset V$, and a function $G : \mathbb{C}^d \rightarrow \mathbb{C}$ such that (1) $G$ is holomorphic on $V$ and (2) $g = G_{|\R^d}$. (This is sometimes known as the {\it complexification} of $g$; see, e.g., Ch. 6.4 of \cite{krantz2002primer}). 

Let $S = B(0, (k+1)r)$. Since $V$ is open and $S$ is compact, there exists $\delta > 0$ such that $B_\delta(S) = \{z : |\text{Im}(z)| < \delta, \text{Re}(z) \in S\} \subset V $.

Define $K : \mathbb{C}^d \times \R^d \rightarrow \mathbb{C}$ as $K(z, y) = G(z-y)\mathbbm{1}_A(y)$; $K$ is holomorphic with respect to $z$ on $B_\delta(B(0, kr))$, since $z-y \in B_\delta(S)$ whenever $y\in B(0, r)$, and $K(z, y) = 0$ otherwise. By Cauchy's integral formula, for any closed curve $\gamma\subset B_\delta(B(0, kr))$, $\int_\gamma K(z, y) dz = 0$. 

Let $F:\mathbb{C}^d \rightarrow \mathbb{C}$ with $F(z) = \int_{\R^d} K(z, y)\, dy$, a function that is continuous in $B_\delta(\R^d)$. Note here that $f = F_{\mid \R^d}$.

Integrating over a closed curve $\gamma\subset B_\delta(B(0, kr))$:
\[\int_\gamma F(z)\, dz = \int_\gamma \int_{\R^d} K(z, y)\, dy\, dz\]
By Fubini's theorem, the order of integration can be exchanged:
\[\int_\gamma F(z)\, dz = \int_{\R^d}\int_\gamma  K(z, y)\,  dz\, dy = 0\]
By Morera's theorem, since the domain $B_\delta(B(0, kr)))$ is simply connected, this is a sufficient condition that $F$ is holomorphic, and thus $f$ is analytic. 
\end{proof}

Since $f$ is analytic, we can argue that $\{x : f(x) = c\}$ has measure 0 in $\R^d$. This holds by the identity theorem for analytic functions; if $f(x) = c$ on an $S\subseteq \R^d$ where $S$ has an accumulation point, then $f =c$ everywhere. Since $g(x)$ vanishes as $\lVert x \rVert\rightarrow \infty$, $f$ cannot be constant. Hence, the level sets have measure 0. 
\end{proof}
In order to apply the results for the general dynamics (Definition~\ref{alg:cts-general}), we need only show that the threshold function fulfills Assumption~\ref{assumption:gamma}. To do this, we must show that $A_t$ has bounded support; that is, there exists an $R$ such that $A_t \subseteq B(0, R)$ for all $t \ge 1$. 

First, we will argue that $A_t$ is bounded by a ball $B(0,r)$, with $r$ dependent on $A_0$, for all $t$. 
\begin{lemma}
\label{lem:bounded0}
    Suppose $|A_0| =\alpha$. Then, there is some $R>0$ such that $A_1\subseteq B(0,R)$.
\end{lemma}
\begin{proof}

    Let $C_0$ be the threshold such that $A_1=\{x\in\R^d:f_0(x)\ge C_0\}$. Choose $\epsilon>0$ such that $\epsilon g(0)<C_0/2$.
    
    Let $r$ be a real number such that $|A_0 \cap B(0, r)|<\epsilon$. Such an $r$ must exist by monotone convergence, as $\lim_{r \rightarrow \infty} \int_{\R^d \setminus B(0, r)} \mathbbm{1}_{A_0}(x)\, dx = 0$. Choose $r'$ so that $g(r')<\dfrac{C_0}{2\alpha}$. Then, for $\lVert x\rVert>r+r'$,
    \begin{align*}
        \int_{\R^d\cap A_0}\mathbbm1{A_0}(y)g(x-y)dy&=\int_{B(0,r)}\mathbbm1{A_0}g(x-y)dy+\int_{\R^d \setminus B(0,r)}\mathbbm1{A_0}g(x-y)dy\\
        &\le \alpha g(r')+\epsilon g(0)\\
        &<C_0
    \end{align*}
 Thus, $A_1\subseteq B(0,r+r')$.
\end{proof}
Next, a universal bound over $A_t$ follows from an earlier lemma, Lemma~\ref{lem:helper0}, which applies to a general class of functions $\psi$. 
\begin{lemma}
    Suppose $A_1\subseteq B(0,R/4)$. Then, $A_t\subseteq B(0,R)$ for all $t\ge 1$.
    \label{lem:bounded}
\end{lemma}
\begin{proof}
    Note that any hyperplane $H$ with $B(0,R/4)\subseteq H^+$ is reflecting, and this will remain true for all $t\ge 1$.
    
    By Lemma~\ref{lem:helper0}, if there is $x\in A_t$ with $\|x\|\ge R$, $B(0,R-R/2)\subseteq A$, and so
    $|A_t|\ge |B(0,R/2)|>|A_1|$, which is impossible.
\end{proof}

With the knowledge that $A_t$ is bounded across time steps, the $\alpha$-cap process (Definition~\ref{alg:cts}) is a special case of the general process (Definition~\ref{alg:cts-general}). Thus, Theorem \ref{thm:general} applies to this process, and implies the following result.
\thmmain*

\section{Discussion}
We have provided a continuous generalization of communication dynamics in a network of binary neurons with a geometric graph structure. Taking the limit as the number of vertices approaches infinity, we model the firing probabilities of individual neurons in $\R^d$ with a sequence of measurable functions $\psi_t:\R^d\rightarrow[0,1]$. The function $\psi_t$ is evolved by convolving with a kernel function $g$, representing the edge probabilities of the graph, followed by sharpening with a monotone function $\gamma_t$, representing the nonlinear gating function used to update the state. We rigorously demonstrated that for a broad set of such processes, the sequence of functions $\psi_t$ tends toward radial symmetry over time. 

As noted in earlier sections, the process studied in this paper is a generalization of the MBO scheme for approximating motion by mean curvature~\cite{merriman1992diffusion}. The similarities in these processes offer a surprising connection between motion by mean curvature and geometric networks of binary neurons. While the two physical systems seem unrelated, they have a conceptual similarity; in both processes, the state of a point in space - either the velocity or the binary state of the vertex - evolves by aggregating the state of neighboring points. Further research into geometric neuron networks could leverage results in interface motion to obtain new insights. 
\section{Conclusion and Future Work}
This work opens up several promising open questions, detailed here.

\paragraph{Extension to general kernels}The results presented in this paper make heavy use of the radial symmetry of the kernel. While Euclidean distance is a common metric used to model spatial networks, it is far from the only possibility. For a variant of the process where $g$ is a function of a different metric --- for example, $L_1$ distance in $\R^d$ --- it is not clear even what the fixed points of the process should look like. 

Another interesting direction is to study the infinite analogue of `hard' geometric random graphs; in this variant, $g(x,y) = 1$ if $\|x-y\|<r$, and 0 otherwise. Our results do not apply to this kernel, since we require that $g$ is nonzero and strictly decreasing in distance. In fact, the process with this kernel does not always converge to radial symmetry. For example, in the fixed-volume $\alpha$-cap process under this kernel, a set of 2 (or more) sufficiently spaced, equally sized balls is a fixed point. Investigating the dynamics here could yield interesting insights. 

\paragraph{Connecting the continuous and discrete models}The process studied in this work is defined as the continuous generalization of a model of binary neurons with a random geometric structure. The dynamics of the two processes are qualitatively similar, both in experiment and in theoretical results~\cite{reid2023k}. A next step may be to determine if the theorems presented here imply high-probability convergence of the $k$-cap process on a broader range of graph models.

\paragraph{Analyzing the rate of convergence}This work implies a lower bound on the rate of convergence; the change in the width of reflecting slabs, $d_t(u) - d_{t+1}(u)$, is lower bounded by a function of $R$, $d$, and parameters dependent on the function $g$ (per Lemma \ref{lem:reflection_progress}). However, the bounds given are not optimal, as they assume that $\psi_t$ has a `worst-case' structure at each time step. There may be other properties of $\psi_t$ which can give a tighter bound on the convergence rate.

\section{Appendix}
This section contains the proof that the set of reflecting slabs have a unique common intersection (Lemma \ref{lem:slab_intersection_unique}).
\begin{lemma}
    \label{lem:slab_intersection}
    All reflecting slabs of $\psi_t$ have a common intersection.
\end{lemma}
\begin{proof}
Let $\{H_\alpha\}_{\alpha\in I}$ be the collection of all reflecting hyperplanes for $\psi_t$. Suppose for the sake of contradiction that the reflecting slabs have an empty intersection. Then, $\bigcap_{\alpha\in I}\overline{H_{\alpha}^+}=\emptyset$.

Let $K$ be the intersection of the reflecting slabs in the directions of the coordinate axes. Note that $K$ is compact. Then,
\begin{align*}
    \bigcap_{\alpha\in I}(\overline{H_{\alpha}^+}\cap K)=\emptyset
\end{align*}
so by the finite intersection property of compact sets, there is a finite subcollection of $\{\overline{H_{\alpha}^+}\cap K\}_{\alpha\in I}$ with empty intersection. Since $K$ itself is the intersection of finitely many hyperplanes $\overline{H_{\alpha}^+}$, there is a finite subcollection of $\{\overline{H_{\alpha}^+}\}_{\alpha\in I}$ with empty intersection. 

    Let $H_i=\{a_i^{\top}x=b_i\}$ for $i=1,\ldots,k$
    be a minimal collection of reflecting hyperplanes such that $\overline{H_1^+}\cap \cdots\cap \overline{H_k^+}=\emptyset$. In other words, $Ax\ge b$ has no solution, where
    \begin{align*}
        A&=\begin{bmatrix}
            a_1^{\top}\\
            \vdots\\
            a_k^{\top}
        \end{bmatrix}&
        b&=\begin{bmatrix}
            b_1\\
            \vdots\\
            b_k
        \end{bmatrix}
    \end{align*}
    but there exists a solution if any constraint is removed. 

    By Farkas' lemma (see e.g.\cite{schrijver1998theory} Corollary 7.1e),
    \begin{align*}
        A^{\top}y&=0\\
        b^{\top}y&=1\\
        y&\ge 0
    \end{align*}
    has a solution $y^*\in\R^k$, but has no solution where $y_i=0$, for any $i$.

    The feasibility of this system implies that $\rank\begin{bmatrix}A&b\end{bmatrix}=k$; or else it would be possible to find a fulfilling $y$ with $y_i = 0$ by perturbing $y^*$. 

    Assume without loss of generality that the first $k-1$ rows of $A$ are linearly independent. By applying an affine transformation, we can assume that
    \begin{align*}
        \overline{H_i^+}&=\{x_i\le 0\}&\text{for $i=1,\ldots,k-1$}\\
        \overline{H_k^+}&=\{c^{\top}x\ge \beta\}
    \end{align*}
    Since $0\notin \bigcap_{j=1}^k\overline{H_j^+}$, $\beta>0$.
    
    We must have $c_i=0$ for $k\le i\le d$, or else the intersection of $\{x_i=0\}$ for $1\le i<k$ and $\{c^{\top}x=\beta\}$ would have a solution.

    For $1 \le i < k$, $c_i\ge 0$; otherwise, for a large enough $M$, we will have $-Me_i\in \bigcap_{j=1}^k\overline{H_j^+}$.
    
    For $1 \le i < k$, if $c_i=0$, then there will still be no solution even if we remove the constraint $x_i\le 0$, contradicting minimality. 
    
    Hence, $c_i>0$ for all $1\le i<k$.

    Performing another affine transformation, we see that we may assume that
    \begin{align*}
        \overline{H_i^+}&=\{x_i\le 0\}&\text{for $i=1,\ldots,k-1$}\\
        \overline{H_k^+}&=\{x_1+\cdots+x_{k-1}\ge k-1+\sqrt{k-1}\}
    \end{align*}
    
    Shifting by $1$ in each coordinate, we may instead assume that
    \begin{align*}
        \overline{H_i^+}&=\{x_i\le -1\}&\text{for $i=1,\ldots,k-1$}\\
        \overline{H_k^+}&=\{\frac{1}{\sqrt{k-1}}(x_1+\cdots+x_{k-1})\ge 1\}
    \end{align*}

    In particular, for all $1\le i\le k$, $\overline{H_i^+}=\{u_i^{\top} x\le -1\}$ for some unit vector $u_i\in\R^d$.

    Also note that $\bigcup_{i=1}^k\{u_i^{\top}x\ge 0\}=\R^d$.

    \begin{claim}
        For any $y\in \R^d$, there is some $i\in\{1,\ldots,k\}$ such that
        \begin{itemize}
            \item $B(y,1)\subseteq H_i^-$, so in particular $\int_{B(y^{H_i},1)} \psi_t(x)\, dx\ge \int_{B(y,1)} \psi_t(x)\, dx$
            \item $\|y^{H_i}\|^2\ge \|y\|^2+4$
        \end{itemize}
    \end{claim}
    \begin{proof}
        For a given $y$, we choose $i$ so that $u_i^{\top}y\ge 0$. Then,
        \begin{align*}
            y^{H_i}&=y-2(u_i^{\top}y+1)u_i\\
            \|y^{H_i}\|^2&=\|y\|^2+\|2(u_i^{\top}y+1)u_i\|^2-4y^{\top}(u_i^{\top}y+1)u_i\\
            &=\|y\|^2+4(u_i^{\top}y+1)^2-4u_i^{\top}y(u_i^{\top}y+1)\\
            &=\|y\|^2+4u_i^{\top}y+4\\
            &\ge \|y\|^2+4
        \end{align*}
    \end{proof}

    Since the support of $\psi_t$ is bounded, $C:=\sup\{\|x\|^2:x\in\R^d,\int_{B(x,1)}\psi_t(y)\, dy>0\}<\infty$.

    Choose $y\in\R^d$ such that $\int_{B(y,1)}\psi_t(x)\, dx>0$ and $\|y\|^2\ge C-1$.

    By the claim, there is some $i$ such that $\int_{B(y^{H_i},1)}\psi_t(x) \, dx\ge \int_{B(y,1)} \psi_t(x)\, dx>0$ and $\|y^{H_i}\|^2\ge \|y\|^2\ge C+3$, which contradicts the definition of $C$ as the supremum over such $\|y\|^2$. This implies the statement of the lemma.
\end{proof}

\begin{proof}[Proof of Lemma~\ref{lem:slab_intersection_unique}]
\leavevmode\\
The intersection of reflecting slabs for $\psi_t$ is compact and nonempty. By the reflection principle, any reflecting slab at $t$ is also reflecting at $t+1$. This implies the intersection of the set of reflecting slabs at $t$ is a nested sequence of compact sets whose diameter (which is at most $d_t$) goes to zero. Such sequences have a unique common point.
\end{proof}
\section*{Acknowledgments}The authors would like to thank Santosh Vempala and Ben Jaye for valuable discussion and feedback. 
M. R. was supported in part by NSF Award CCF-2106444 and a Simons Investigator award. D. Z. was supported by NSF Award CCF-2338816.
\printbibliography
\end{document}